\documentclass[aps,onecolumn,10pt,a4paper,nofootinbib,floatfix]{revtex4}
\usepackage{graphicx,amsmath}
\usepackage{amssymb}
\usepackage{amsthm}
\usepackage{hyperref}
\usepackage{doi}
\usepackage{centernot}
\usepackage{makecell}
\usepackage{multirow}
\usepackage{xcolor}
\usepackage[export]{adjustbox}
\allowdisplaybreaks[4]

\newcommand{\fpdim}[1]{{\rm FPdim}(#1)}

\begin{document}
\title{On the unitarity and modularity of ribbon tensor categories associated with affine Lie algebras}
\author{Daria Rudneva}
\author{Eddy Ardonne}
\affiliation{Department of Physics, Stockholm University, SE-106 91 Stockholm, Sweden}

\begin{abstract}
We study the unitarity and modularity of ribbon tensor categories derived from simple affine
Lie algebras, via their associated quantum groups. Based on numerical calculations,
and assuming two conjectures, we provide the complete picture for which values of $q$
these ribbon tensor categories are (pseudo-)unitary and for which values of $q$ they are
modular. We compare our results with the extensive rigorous results appearing in the literature,
finding complete agreement. For the cases that do not appear in the literature, we complete
the picture.
\end{abstract}
\maketitle

\section{Introduction}
\label{sec:introduction}
Non-abelian anyons were first theorized to exist in (two-dimensional) fractional quantum Hall systems \cite{moore91}.
Non-abelian particles can occur in two-dimensional topological phases, whose low-energy
effective theories are described by topological quantum field theories \cite{witten88}. 
The structure underlying these topological quantum field theories are modular tensor categories \cite{turaev92}.
Later on, it was realized that anyons can be used for fault tolerant quantum computation \cite{kitaev03}
(see also \cite{freedman03,nayak08}).

It is clearly interesting to determine which types of non-abelian particles can exist in principle in
two-dimensional systems. Answering this questions fully is out of reach, even though classification results
exists when the number of particles is low \cite{rowell09}.
One avenue to construct modular tensor categories is by means of quantum groups associated
with affine Lie algebras, see f.i. \cite{rowell06}. In this way, one obtains so-called ribbon tensor categories (RTCs).
In order that one obtains a system describing anyons that could in principle be used in the context of
quantum computation, these categories should be both modular and unitary.

Many (rigorous) results on the modularity and unitarity of RTCs have appeared in the literature \cite{rowell06}, but 
the picture is not complete. In this paper, we provide the full picture, based on numerical calculations and
under the assumption of two conjectures. Our results are fully consistent with the rigorous results that appeared in
the literature. In the remaining cases, our results are clearly conjectural.

In Sec.~\ref{sec:background}, we provide some minimal background, setting the notation. In Sec.~\ref{sec:uni-mod},
we explain how we determine which RTCs are unitary, and which are modular and which conjectures we assume.
The results are presented in Sec.~\ref{sec:results}, which are compared to the existing literature. We close
with a brief discussion in Sec.~\ref{sec:background}, while in the appendices we give some basic Lie algebra
data, and indicate in detail which cases we checked numerically.

\section{Background}
\label{sec:background}

In this section, we give a very brief introduction to (affine) Lie algebras, their associated quantum groups
and how to obtain the categorical data, mainly to set the notation.
This material is standard, we refer to references \cite{book:kassel, book:fuchs, book:gomez, book:bakalov, kitaev06, book:wang, BYB} for more detail.
For the construction of ribbon tensor categories (RTC) from Lie algebras, we refer to \cite{rowell06} in particular.
More information on tensor categories in general can be found in \cite{book:bakalov}.

Let $\mathfrak{g}_r$ denote a Lie algebra of type $A$ - $G$ with rank $r$
(as usual, $A$ - $D$ are the four infinite series and $E$, $F$ and $G$ the exceptional cases).
We denote the simple roots and fundamental weights by $\alpha_i$ and $\omega_i$ respectively, where $1\leq i \leq r$.
The inner product between weights is denoted by $\langle , \rangle$ and given in terms
of the quadratic form matrix, namely $\langle \omega_i , \omega_j \rangle = F_{i,j}$.
Here $F_{i,j}$ are the elements of the quadratic form matrix, which can be written as
$F_{i,j} = (A^{-1})_{i,j}\frac{\alpha_j^2}{2}$, where $A$ denotes the Cartan matrix of $\mathfrak{g}_r$.
In particular, $\langle , \rangle$ is normalized such that $\langle \alpha, \alpha \rangle = 2$ when $\alpha$ is a {\em long} root%
\footnote{We note that Rowell \cite{rowell06} uses a different normalization, namely such that $\langle \alpha, \alpha \rangle = 2$
when $\alpha$ is a {\em short} root.}.
As usual, $\rho = (1,1,\ldots, 1)$ is half the sum of the positive roots.
Finally, the $q$-deformed numbers are defined via $[n]_q = \frac{q^n - q^{-n}}{q - q^{-1}}$. We collect some basic Lie 
algebra data in Appendix~\ref{app:Lie-algebra-data}.

\subsection{From quantum groups to $F$ and $R$ symbols}

The quantum groups associated with affine Lie algebras, such that $q$ is a primitive root of unity, give rise to
ribbon tensor categories (RTCs), which are also called premodular categories.
That is, one obtains $F$-symbols that satisfy the pentagon equations, and $R$-symbols that satisfy the
hexagon equations. In addition, one also obtains a pivotal structure consistent with the braiding. Thus, one obtains
a RTC. We note that the pivotal structure is in fact is spherical.
The RTC one obtains may or may not be modular. If it is modular, one in fact has a modular tensor category (MTC).
The RTC may or may not be unitary, independent whether or not it is modular.
If the RTC is both unitary and modular, one calls the tensor category a unitary modular tensor category (UMTC).

To explicitly construct a RTC from a quantum group, one has to calculate the categorical data.
On starts by specifying the deformation parameter $q$, which is a primitive root of unity. 
We parametrize $q$ by two parameters
$\ell$ and $p$ as follows, $q = e^{\frac{2\pi i p}{\ell}}$, where $\gcd(p,\ell) = 1$.
The allowed values of $\ell$ depend on the type of algebra, and its rank.
Here, the ratio between the length of the long and short roots, denoted by $m$, plays an important role.
For the algebras of type $A$, $D$ and $E$, we have $m=1$. For types $B$, $C$ and $F$ we have $m=2$ while for
type $G$, we have $m=3$.
The RTCs for which $m \mid \ell$ are called uniform, while those for which $m \centernot \mid \ell$ are called
non-uniform.

For the uniform cases, it is customary (at least in the physics literature)
to use the level $k$ instead of $\ell$ to parametrize $q$. In the tables in Sec.~\ref{sec:results} where
we summarize our results, we give both labels. The relation between $\ell$ and $k$ is given by
\begin{equation}
\ell = m (k + g) \ ,
\end{equation}
where $g$ is the dual coxeter number of the algebra.
Because $m \mid \ell$ for the uniform cases and $m \centernot \mid \ell$ for the non-uniform cases,
we have that $k$ is integer for the uniform cases, and non-integer for the non-uniform cases.

Given the type and rank of the algebra, as well as the root of unity (in particular, the value of $\ell$), we need to
specify the simple objects. To do so, we first define the set of dominant weights, denoted by $P_+$ as the
set of weights that are a non-negative integer combination of the fundamental weights.
For the uniform cases $m \mid \ell$, the set of simple objects $C$ is given by
\begin{equation}
C = \{ \lambda \in P_+ | \langle \lambda + \rho, \vartheta \rangle < \ell \} \ ,
\end{equation}
where $\vartheta$ is the highest root. For the non-uniform cases 
$m \mid \ell$, the set of simple objects $C$ is given by
\begin{equation}
C = \{ \lambda \in P_+ | \langle \lambda + \rho, \vartheta_s \rangle < \ell \} \ ,
\end{equation}
where $\vartheta_s$ is the highest {\em short} root. Because the simple objects are associated with
dominant weights, or irreducible representations, we also refer to simple objects as `irreps'.

In the context of the Lie algebra, one can decompose the tensor product of two irreducible representations as a sum of irreducible representations.
In the context of the quantum group with $q$ a primitive root of unity, the tensor product is replaced by the fusion product,
\begin{equation}
\label{eq:fusion-rules}
\lambda \times \mu = \sum_{\nu \in C} N_{\lambda,\mu}^{\nu} \nu \ ,
\end{equation}
where the fusion coefficients $N_{\lambda,\mu}^{\nu}$ are non-negative integers. We denote the trivial simple object by $\mathbf{1} = (0,0,\ldots,0)$ and the (unique) simple object dual to $\lambda$ by $\lambda^*$. We then have
\begin{align}
N_{\lambda,\mu}^{\nu} &=
N_{\mu,\lambda}^{\nu} =
N_{\lambda,\nu^*}^{\mu^*} =
N_{\lambda^*,\mu^*}^{\nu^*}
&
N_{\lambda,\mu}^{\mathbf{1}} &= \delta_{\mu,\lambda^*} \ .
\end{align}
To each simple object, we associate a fusion matrix $N_\lambda$ whose elements are $(N_\lambda)_{\mu,\nu} = N_{\lambda,\mu}^{\nu}$.
The fusion rules are associative, $N_\lambda N_\mu = N_\mu N_\lambda$, which allows one to define the $F$-symbols.

With the set of simple objects and the fusion rules in place, one can calculate the $F$- and $R$-symbols, as well as the pivotal structure.
We only very briefly outline this procedure here, and refer to the paper \cite{as10}, where the procedure
to obtain the $F$- and $R$-symbols is described in detail.
We note that the methods described in \cite{as10} have been implemented in Mathematica \cite{alatc}.

One starts out by calculating the ($q$-deformed) Clebsch-Gordan (qCG) coefficients.
To this end, one decomposes the tensor product of two irreps in terms of all the possible
irreps, by making use of the co-product. The qCG coefficients can in turn be used to calculate
the $F$-symbols, which are simply the $q$-deformed 6-j symbols obtained from the qCG coefficients.

From the qCG coefficients, one can also calculate the $R$-symbols, using the method described in
\cite{as10}. Alternatively, one can try to directly solve the hexagon equations using the
calculated $F$-symbols, which is relatively straightforward.

Finally, one can obtain the pivotal (or rather, spherical) structure, because as we explain in Sec.~\ref{sec:unitarity}, one can obtain the categorical dimensions $d(\lambda)$ directly from the Lie algebra data and the pivotal coefficients $\varepsilon_\lambda$ are determined by
$
\varepsilon_\lambda = d(\lambda) F^{\lambda,\lambda^*,\lambda}_{\lambda;\mathbf{1},\mathbf{1}}
$.

\section{On unitarity and modularity}
\label{sec:uni-mod}

The main focus of the paper will be on the unitarity and modularity of the RTCs derived from
quantum groups associated with affine Lie algebras.
In both cases, we will assume the validity of a particular conjecture (described below), which greatly simplifies the
calculations. We start by discussing unitarity in coming subsection, followed by a discussion on modularity.
In Sec.~\ref{sec:results}, we present the result, and discuss the connection with the literature, because many
(rigorous) results have been obtained already. 

\subsection{Unitarity}
\label{sec:unitarity}

To determine if a RTC is unitary, it suffices to show that there is a gauge such that
the $F$-matrices are unitary. Obtaining the $F$-matrices is, however, a numerically costly procedure,
even for the RTCs associated with affine Lie algebras, where one can use the underlying quantum group
structure.
We therefore follow a different route, exploiting the fact that it is relatively straightforward to obtain
the categorical (or quantum) dimensions of the simple objects.

We start the discussion by introducing the Frobenius-Perron
dimension associated with an irrep, $\fpdim{\lambda}$, which is the largest (Frobenius-Perron) eigenvalue
of the fusion matrix $N_\lambda$. In particular, $\fpdim{\lambda}$ is positive.
One defines the Frobenius-Perron dimension of the category, $\fpdim{\mathcal{C}}$ to be the sum of the squares of the
Frobenius-Perron dimensions of the irreducible representations,
\begin{equation}
    \fpdim{\mathcal{C}} = \sum_{\lambda\in C} \fpdim{\lambda}^2 \ .
\end{equation}
Apart from the Frobenius-Perron dimension, one also defines the categorical dimension, or the
quantum dimension of an irreducible representation, $d(\lambda)$ as the value of a loop with label $\lambda$.
The quantum dimension of the category $d(\mathcal{C})$ is again the sum of the squares of the quantum
dimensions of the irreducible representations,
\begin{equation}
    d(\mathcal{C}) = D^2 = \sum_{\lambda\in C} d(\lambda)^2 \ .
\end{equation}

In the case of RTCs derived from quantum groups associated with affine Lie algebras, there is a rather
explicit formula to calculate the quantum dimensions for a given parameter $q$.
In particular, we have
\begin{equation}
\label{eq:qdims}
d(\lambda) = \prod_{\alpha \in \Delta^+} \frac{[(\lambda+\rho, \alpha)]_q}{[(\rho, \alpha)]_q} \ ,
\end{equation}
where the product is over all the positive roots.
Given a Lie algebra and the value of $q$, one can in principle calculate the quantum dimensions straightforwardly.
This is useful for the following reason.

In a unitary fusion category, all the quantum dimensions are positive.
It has been conjectured, see for instance \cite{book:wang}, that the opposite is also true,
namely if all the categorical dimensions are positive, then the fusion category is unitary.
Because it is in principle straightforward to compute the quantum dimensions, one can determine if
the category is unitary, under the assumption that this conjecture holds.
In \cite{rowell08}, this procedure was formalized by introducing the concept of {\em pseudo-unitarity}.
A category $\mathcal{C}$ is called pseudo-unitary if $d(\mathcal{C}) = \fpdim{\mathcal{C}}$.
In this case, one has that $d(\lambda) = \pm \fpdim{\lambda}$ for all labels $\lambda$. Moreover, there
then also exists a spherical structure such that $d(\lambda) > 0$ for all $\lambda$.   

To summarize, our strategy to determine whether or not a RTC associated with an affine Lie algebra
is thus as follows. We assume the conjecture that if the categorical dimensions are positive, the category
is unitary. We calculate the categorical dimensions using Eq.~\eqref{eq:qdims}, which is relatively
straightforward.
To determine if a RTC is (pseudo-)unitary, ordinary machine precision calculations of the categorical dimension suffices.
From the categorical dimensions, it is straightforward to determine if the category is pseudo-unitrary or even unitary. If the category is only pseudo-unitary, there exists a spherical structure, such that the category
is unitary.

In the next subsection, we will need the categorical dimensions to higher precision, in order to establish modularity.
For that purpose, we use Mathematica with arbitrary precision calculations with 30 digits
precision to avoid stability issues.

\subsection{Modularity}
For the uniform RTCs associated with affine Lie algebras, one can determine the $S$-matrix in the case
that $q$ is the primitive root of unity with $p=1$, because it has been shown that these are always modular (see below). 
The $S$-matrix can be written in terms of a sum over the Weyl group
of the Lie algebra. With the $S$-matrix at hand, one can determine the fusion rules using the
Verlinde formula \cite{verlinde88}. Here, we assume that the $S$-matrix
is normalized such that $S_{\mathbf{1},\mathbf{1}} = 1$,
\begin{equation}
\label{eq:verlinde}
N_{\lambda,\mu}^{\nu} = \sum_{\kappa \in C} \frac{S_{\lambda,\kappa}S_{\mu,\kappa}S^*_{\nu,\kappa}}{D^2 S_{\mathbf{1},\kappa}} \ .
\end{equation}

Using the fusion rules, one can now even obtain the $S$-matrix for arbitrary roots of unity,
via~\cite{book:bakalov}
\begin{equation}
    S_{\lambda,\mu} = \frac{1}{\theta_\lambda \theta_\mu} \sum_{\nu\in C} N^\nu_{\lambda^*,\mu} d(\nu) \theta_\nu \ ,
\end{equation}
where $\theta_\lambda$ are the twist factors, given by $\theta_\lambda = q^{\frac{\langle\lambda,\lambda+\rho\rangle}{2}}$
for the RTCs we consider. Given the $S$-matrix, it is of course straightforward
to determine if it is invertible, and hence if the RTC is modular.

There is, however, a problem with this procedure. To determine the initial $S$-matrix,
one has, as we indicated, to perform a sum over the Weyl group. This severely restricts the cases
one can actually calculate numerically%
\footnote{We note that using the quantum Racah formula \cite{sawin06} reduces the complexity.}.
Because of this reason, we use a different strategy. Following \cite{book:wang}, we define
the quantities
\begin{equation}
\label{eq:p-relation}
    p_{\pm} = \sum_{\lambda\in C} \theta^{\pm1}_\lambda d(\lambda)^2 \ .
\end{equation}
Then, if the category is modular, one has $p_+ p_- = D^2 = \sum_{\lambda\in C} d(\lambda)^2$ \cite{book:wang}.
Using the program \cite{alatc}, we observed that for all cases checked, the converse is also true. Namely,
having the relation $p_+ p_- = D^2$ is sufficient for the RTC to be modular.

This motivates the following definition. We call a RTC for which the relation
$p_+ p_- = D^2$ is satisfied {\em pseudo-modular}. We conjecture that for RTCs
associated with simple affine Lie algebras, pseudo-modular categories are modular.

Under the assumption that our conjecture holds, it is relatively straightforward to determine if a RTC associated with an affine Lie algebra is modular. In the previous subsection, we explained
how to obtain the categorical dimensions $d(\lambda)$. It is therefore straightforward to determine if the
relation Eq.~\eqref{eq:p-relation} is satisfied. Note that because of the sum over the simple objects in the definition of
$p_\pm$, it is necessary to use arbitrary precision calculations to avoid stability issues.
We used Mathematica to perform these calculations, using
30 digits precision, which suffices in all cases we considered.
 
\section{Results}
\label{sec:results}
In this section, we give the conditions such that the RTCs associated with
affine Lie algebras are unitary and or modular. We first consider the uniform cases,
followed by the non-uniform cases.

We recall that the primitive root of unity $q$ is parametrized by the two parameters
$\ell$ and $p$ as $q = e^{\frac{2\pi i p}{\ell}}$.
The uniform cases correspond to $m \mid \ell$, while for the non-uniform cases, $m \centernot \mid \ell$.
Because of the relation $\ell = m (k + g)$, we have that the level $k$ is integer for the uniform cases,
and non-integer for the non-uniform cases.

To have an admissible value of $q$, the condition $\gcd(p,\ell) = 1$ needs to be fulfilled.
In the tables below, we do not explicitly state this condition every time, but it is implicitly assumed to be valid.
Sending $q \rightarrow q^{-1}$ does not change the unitarity and modularity of the RTC.
We therefore choose $p$ to lie in the range $-\ell/2 < p < \ell/2$, and we only report the results for
positive $p$, in the range $1 \leq p < \ell/2$.

Below, we compare our results with the (incomplete) results that appear in the literature.
In particular, we compare with the results that appear in the paper by Rowell \cite{rowell06}.
We note that the definition of $q$ in \cite{rowell06} (denoted here by $q_{\rm Rowell}$)
differs from the definition that we use in the current paper (denoted by $q = q_{\rm ours}$).
The relation is given by $q_{\rm ours} = q_{\rm Rowell}^2$.
Moreover in \cite{rowell06}, $q_{\rm Rowell}$ is parametrized as $q_{\rm Rowell} = e^{\frac{\pi i z}{\ell}}$, which
means that $p = z$. The reason to use a different symbol, namely $p$, appearing in the parametrization of
$q$ is to remind the reader that we use a different definition of $q$.

\subsection{Uniform cases}
We give the results for the unitarity and modularity of the RTCs, for the
{\em uniform cases} in Table~\ref{tab:uni-mod-uniform}.
We indicate if the theories are unitary, pseudo-unitary and modular (under the assumption that the
conjectures mentioned in sec.~\ref{sec:uni-mod} hold). In the cases where the category is pseudo-unitary
but not unitary, the pseudo-unitarity is indicated in red. In these cases, there exists a different
spherical structure such that all categorical dimensions are positive. For that spherical structure, the
category is also unitary.
It was shown that for $p=1$, the 
RTCs are both unitary and modular. Our results are consistent with these results.
In subsection~\ref{sec:literature-comparison} below, we compare our results with the existing literature in more detail.
We note that $(a | n)$, where $a$ is an integer, and $n$ a non-negative odd integer,
denotes the Jacobi symbol. We comment on the appearance of the Jacobi symbols in
sec.~\ref{sec:literature-comparison} below.

\begin{table}[ht]
\begin{tabular}{|c|c|c|c|c|c|}
\hline
& $k$ & $\ell$ & unitary & pseudo-unitary & modular \\
\hline\hline
\multirow{2}{*}{$A_{r}$} & $k = 1$ & $\ell = r + 2$ & $2\centernot \mid p$ & \color{red}$\forall p$ & \multirow{2}{*}{$\gcd(p,r+1)=1$} \\
\cline{2-5}
& $k \geq 2$ & $\ell \geq r + 3$ & $p=1$ & $p=1$ & \\
\hline\hline
\multirow{3}{*}{$B_{r}$} & $k=1$ & $\ell = 4r$ & $p = 1, 7 \bmod 8$ & \color{red}$\forall p$ & \multirow{3}{*}{$\forall p$} \\
\cline{2-5}
& $k=2$ & $\ell = 4r+2$ & $( 2r+1 | p ) = 1$ & \color{red}$\forall p$ & \\
\cline{2-5}
& $k\geq 3$ & $\ell \geq 4r+4$ & $p=1$ & $p=1$ & \\
\hline\hline
\multirow{3}{*}{$C_{r}$} & $k = 1$ & $\ell = 2r+4 $ & \multirow{3}{*}{$p=1$} &
\makecell{$p = 1$ \\ \color{red}$2 \mid r \wedge p = r + 1$}
& \multirow{3}{*}{$\forall p$} \\
\cline{2-3} \cline{5-5}
& $k = 2$ & $\ell = 2r + 6$ & & 
\makecell{$p = 1$ \\ \color{red}$r = 2 \wedge p = 3$} & \\
\cline{2-3} \cline{5-5}
& $k \geq 3$ & $\ell \geq 2r+8$ & & $p=1$ & \\
\hline\hline
\multirow{3}{*}{$D_{r}$} & $k=1$ & $\ell = 2r-1$ &
\makecell{$r = 0, 1 \bmod 4 \wedge \forall p $\\ $r = 2, 3\bmod 4 \wedge 2\centernot \mid p$}
& \color{red}$\forall p$
& \multirow{3}{*}{$2 \centernot \mid p$} \\
\cline{2-5}
& $k=2$ & $\ell = 2r$ & $(r|p) = 1$ & \color{red}$\forall p$ & \\
\cline{2-5}
& $k\geq 3$ & $\ell \geq 2r + 1$& $p=1$ & $p=1$ & \\
\hline\hline
\multirow{4}{*}{$E_{6}$} & $k=1$ & $\ell = 13$ & $\forall p$ & $\forall p$ & \multirow{4}{*}{$3\centernot \mid p$} \\
\cline{2-5}
& $k=2$ & $\ell = 14$ & $p=1$ & $p=1$ & \\
\cline{2-5}
& $k=3$ & $\ell = 15$ & $p=1, 4$ & $p=1, 4$ & \\
\cline{2-5}
& $k\geq 4$ & $\ell \geq 16$ & $p=1$ & $p=1$ & \\
\hline
\multirow{4}{*}{$E_{7}$} & $k=1$ & $\ell = 19$ & $2\centernot \mid p$ & \color{red}$\forall p$ & \multirow{4}{*}{$2 \centernot \mid p $} \\
\cline{2-5}
& $k=2$ & $\ell = 20$ & $p=1, 9$ & $p=1, 9$ & \\
\cline{2-5}
& $k=3$ & $\ell = 21$ & $p=1, 5$ & $p=1, {\color{red}4}, 5$ & \\
\cline{2-5}
& $k\geq 4$ & $\ell \geq 22$ & $p=1$ & $p=1$ & \\
\hline
\multirow{6}{*}{$E_{8}$} & $k=1$ & $\ell = 31$ & $\forall p$ & $\forall p$ & \multirow{6}{*}{$\forall p$} \\
\cline{2-5}
& $k=2$ & $\ell = 32$ & $p=1, 7,  9, 15$ & \color{red}$\forall p$ & \\
\cline{2-5}
& $k=3$ & $\ell = 33$ & $p=1, 10$ & $p=1, 10$ & \\
\cline{2-5}
& $k=4$ & $\ell = 34$ & $p=1$& $p=1$& \\
\cline{2-5}
& $k=5$ & $\ell = 35$ & $p=1, 6 $& $p=1, 6 $& \\
\cline{2-5}
& $k\geq 6$ & $\ell \geq 36$ & $p=1$& $p=1$& \\
\hline\hline
\multirow{5}{*}{$F_{4}$} & $k=1$ & $\ell = 20$ & $p=1, 9$ & $p=1, 9$ & \multirow{5}{*}{$\forall p$} \\
\cline{2-5}
& $k=2$ & $\ell = 22$ & $p=1$ & $p=1$ & \\
\cline{2-5}
& $k=3$ & $\ell = 24$ & $p=1, 5$ & $p=1, 5$ & \\
\cline{2-5}
& $k=4$ & $\ell = 26$ & $p=1, 5$ & $p=1, 5$ & \\
\cline{2-5}
& $k\geq 5$ & $\ell \geq 28$ & $p=1$ & $p=1$ & \\
\hline\hline
\multirow{5}{*}{$G_{2}$} & $k=1$ & $\ell = 15$ & $p=1,4$ & $p=1,4$ & \multirow{5}{*}{$\forall p$} \\
\cline{2-5}
& $k=2$ & $\ell = 18$ & $p=1$ & $p=1$ & \\
\cline{2-5}
& $k=3$ & $\ell = 21$ & $p=1, 4, 5$ & $p=1, 4, 5$ & \\
\cline{2-5}
& $k=4$ & $\ell = 24$ & $p=1, 5$ & $p=1, 5$ & \\
\cline{2-5}
& $k\geq 5$ & $\ell \geq 27$ & $p=1$ & $p=1$ & \\
\hline
\end{tabular}
\caption{The (pseudo-)unitarity and modularity conditions for the uniform cases. The cases which are pseudo-unititary but not unitary are indicated in red. We note that the condition $\gcd(p,\ell)=1$ is always implicitly assumed.}
\label{tab:uni-mod-uniform}
\end{table}

\subsection{Non-uniform cases}

We give the results for the unitarity and modularity of the RTCs, for the
{\em non-uniform cases} in Table~\ref{tab:uni-mod-non-uniform}. The notation is the same as for
the uniform cases in the previous subsection.

\begin{table}[ht]
\begin{tabular}{|c|c|c|c|c|c|}
\hline
& $k$ & $\ell$ & unitary & pseudo-unitary & modular \\
\hline\hline
\multirow{3}{*}{$B_{r}$} & $k = (-2r+3)/2$ & $\ell = 2r+1$ & 
\makecell{
$r = 0 \bmod 4 \wedge \forall p$ \\
$r = 1 \bmod 4 \wedge 2 \centernot \mid p$ \\
$r = 3 \bmod 4 \wedge 2 \mid p$
}
& \color{red}$\forall p$
& \multirow{3}{*}{$2\centernot \mid r \wedge 2\centernot \mid p $}  \\
\cline{2-5}
& $k = (-2r+5)/2$ & $\ell = 2r+3$ &
\makecell{
$2 \centernot \mid r \wedge p = (r+1)/2$ \\ $4 \mid r \wedge p = (r+2)/2$
}
&
\makecell{
$2 \centernot \mid r \wedge p = (r+1)/2$ \\ \color{red}$2 \mid r \wedge p = (r+2)/2$
}
& \\
\cline{2-5}
& $k \geq (-2r+7)/2$ & $\ell \geq 2r+5$ & $-$ & $-$ & \\
\hline\hline
\multirow{3}{*}{$C_{r}$} & $k = -1/2$ & $\ell = 2r+1$ & \multirow{3}{*}{$-$} & \color{red}$\forall p$ & \multirow{3}{*}{$-$} \\
\cline{2-3} \cline{5-5}
& $k = 1/2$ & $\ell = 2r+3$ & & \color{red}$p=2$ & \\
\cline{2-3} \cline{5-5}
& $k \geq 3/2$ & $\ell \geq 2r+5$ & & $-$ & \\
\hline\hline
\multirow{4}{*}{$F_{4}$} & $k = -5/2$ & $\ell = 13$ & $\forall p$ & $\forall p$ & \multirow{4}{*}{$\forall p$ }  \\
\cline{2-5}
&  $k = -3/2$ & $\ell = 15$ & $-$ & $-$ & \\
\cline{2-5}
&  $k = -1/2$ & $\ell = 17$ & $p=3$ & $p=3$ & \\
\cline{2-5}
&  $k \geq 1/2$ & $\ell \geq 19$ & $-$ & $-$ & \\
\hline\hline
\multirow{7}{*}{$G_{2}$} & $k = -5/3$ & $\ell = 7$ & $\forall p$ & $\forall p$ & \multirow{7}{*}{$\forall p$} \\
\cline{2-5}
& $k = -4/3$ & $\ell = 8$ & $-$ & \color{red}$\forall p$ & \\
\cline{2-5}
& $k = -2/3$ & $\ell = 10$ & $-$ & $-$ & \\
\cline{2-5}
& $k = -1/3$ & $\ell = 11$ & $p=2$ & $p=2$ & \\
\cline{2-5}
& $k = 1/3$ & $\ell = 13$ & $p=3$ & $p=3$ & \\
\cline{2-5}
& $k = 2/3$ & $\ell = 14$ & $p=3$ & $p=3$ & \\
\cline{2-5}
& $k \geq 4/3$ & $\ell \geq 16$ & $-$ & $-$ & \\
\hline
\end{tabular}
\caption{The (pseudo-)unitarity and modularity conditions for the non-uniform cases. The cases which are pseudo-unititary but not unitary are indicated in red. We note that the condition $\gcd(p,\ell)=1$ is always implicitly assumed.}
\label{tab:uni-mod-non-uniform}
\end{table}

\subsection{Comparison with the literature}
\label{sec:literature-comparison}

In this section, we compare the results for the uniform and non-uniform cases, as given in the previous two subsections to results that appeared in the literature.
Perhaps the most important result is that for all uniform cases with $p=1$, the
RTC is both unitary \cite{wenzl98,xu98} and modular \cite{kirillov96}.
We mainly compare our results on modularity with those given in \cite{rowell06},
and on unitarity with those given in \cite{rowell08}.
As usual, the condition $\gcd(p,l)=1$ is assumed implicitly.

\subsubsection{Modularity of the uniform cases}
The modularity of the uniform cases that were completely covered
(that is, for all possible values of $q$) in \cite{rowell06} are
$A_r$, $B_r$, $C_r$, $E_8$, $F_4$ and $G_2$. Our
results agree with the results for these cases.

In addition, the cases $D_r$ with $l$ even, $E_6$ with $l$ a multiple of
$3$ and $E_7$ with $l$ even were covered for all possible values of $p$.
Again, our results agree with these results.

Finally, the cases $D_r$ with $l$ odd and $p$ odd, $E_6$ with $l$ not a multiple of $3$
and $p$ not a multiple of $3$ and $E_7$ with $l$ odd and $p$ odd were shown to be
modular. Again, our result are in agreement. For the remaining values of $p$, it is not
known if the RTCs are modular, but \cite{rowell06} expressed the expectation that these
cases are not modular. Our results are in agreement with that expectation.

\subsubsection{Modularity of the non-uniform cases}
The modularity of the non-uniform cases that were completely covered
in \cite{rowell06} are $B_r$ and $C_r$. Our results agree with those results.
The modularity of the non-uniform cases $F_4$ and $G_2$ is open,
but \cite{rowell06} expects these to be `sometimes' modular. Our results indicate
that they are always modular.

\subsubsection{Unitarity of the uniform cases}
The study of the unitarity of the RTCs associated with
affine Lie algebras seem to have predominantly been concentrated on the non-uniform
cases, see below.

There is the generic case $p=1$, which is always unitary, as discussed above.
Generically, for large enough level $k$ (or large enough $\ell$), $p=1$ is the only unitary case.
The $E$ cases are always pseudo-unitary for $k=1$ (their fusion rules are
equivalent to $\mathbb{Z}_3$, $\mathbb{Z}_2$, and $\mathbb{Z}_1$, for $E_6$, $E_7$ and $E_8$
respectively), while $E_7$ is also pseudo-unitary for $k=2$ (and equivalent to $so(3)_3 \times su(2)_2$).
Apart from these cases, there are a few (pseudo)-unitary cases for the
exceptional Lie algebras at low level.

Perhaps the most interesting cases are $B_r$ and $D_r$ at level $k=1$ and in particular,
level $k=2$. The $B_r$ case for $k=2$ were studied in detail in \cite{ardonne16,ardonne21}, while
$D_r$ at level $k=2$ was studied in detail in \cite{bruillard19}. In particular, in
\cite{ardonne21}, explicit expressions for the $F$ and $R$-symbols were obtained.
Interestingly, some of the $R$-symbols depended on Jacobi symbols. One reason for this
is lies in the form of a hexagon equation, which takes the form of a quadratic Gauss sum,
leading to Jacobi symbols. These also appear in the the condition on $p$ and $r$, determining
which RTCs of type $B_r$ and $D_r$ at level $k=2$ are unitary
(though they are always pseudo-unitary).

\subsubsection{Unitarity of the non-uniform cases}
In \cite{rowell08}, it is shown that the non-uniform cases of $B_r$ and $C_r$ are
not (pseudo-)unitary for $l \geq 2r+5$. Our results are in agreement with these bounds.
We note that the results of \cite{rowell08} are actually (much) stronger, there exist no
RTC with fusion rules equivalent to those of the non-uniform cases with $l \geq 2r+5$ that
are unitary. 

The cases $F_4$ and $G_2$ are also covered in \cite{rowell08}. For $F_4$, no cases with $l \geq 19$ are
(pseudo-)unitary, while the case $l=17$, $p=3$ is, in agreement with our results (note
that the case $l=13$ is completely trivial).

For $G_2$, no cases with $l\geq 16$ are (pseudo-)unitary, but the cases
$(l,p) = (11,2),(13,3),(14,3)$ are, in agreement with our results (note that the case $l=7$ is completely trivial).

\section{Discussion}
\label{sec:discussion}

Based on extensive numerical calculations, and assuming two conjectures,
we determined which RTCs associated with
affine Lie algebras are (pseudo-)unitary and/or modular. We conjecture that
the results presented here indeed cover all the unitary and modular cases.

Many cases were studied in the literature before, leading to rigorous results.
Our calculations are consistent with the results that previously appeared in the
literature. In the remaining cases, we completed the picture, despite the fact that
we do not provide a proof for these cases.

{\em Acknowledgments ---}
EA would like to thank Eric Rowell, Steve Simon, Joost Slingerland, Gert Vercleyen and Zhenghan Wang for numerous illuminating discussions.

\appendix

\section{Some Lie algebra data}
\label{app:Lie-algebra-data}

In this appendix, we collect some Lie algebra data that is relevant for the paper.
We use the following constraints on the ranks for the four infinite series:
$A_{r\geq 1}$, $B_{r\geq 3}$, $C_{r\geq 2}$ and $D_{r\geq 4}$.
In table~\ref{tab:lie-algebra-data}, we give the dimension, dual Coxeter number $g$,
the highest root $\vartheta$ and the highest short root $\vartheta_s$ in the case of the
non simply-laced Lie algebras. We note that for $A_1$, the highest root is given by 
$\vartheta = (2)$.
We also give $m$, which denotes the ratio of the square
length of the long over the short roots, and (in the non-simply laced cases $B$, $C$,
$F$ and $G$), the short simple roots.

\begin{table}[ht]
\begin{tabular}{|c|c|c|c|c|c|c|}
\hline
type & dimension & $g$ & $\vartheta$ & $\vartheta_s$ & $m$ & short simple roots \\
\hline
$A_{r\geq 1}$ & $r^2+2r$ & $r+1$ & $(1,0,\ldots, 0, 1)$ & - & 1 & - \\
\hline
$B_{r\geq 3}$ & $2r^2+r$ & $2r-1$ & $(0,1,0,\ldots,0)$ & $(1,0,0,\ldots,0)$ & 2 & $\alpha_r$\\
\hline
$C_{r\geq 2}$ & $2r^2+r$ & $r+1$ & $(2,0,\ldots,0)$ & $(0,1,0,\ldots,0)$ & 2 & $\alpha_1,\ldots, \alpha_{r-1}$\\
\hline
$D_{r\geq 4}$ & $2r^2-r$ & $2r-2$ & $(0,1,0,\ldots,0)$ & - & 1 & - \\
\hline
$E_{6}$ & 78 & $12$ & $(0,0,0,0,0,1)$ & - & 1 & - \\
\hline
$E_{7}$ & 133 & $18$ & $(1,0,0,0,0,0,0)$ & - & 1 & - \\
\hline
$E_{8}$ & 248 & $30$ & $(1,0,0,0,0,0,0,0)$ & - & 1 & - \\
\hline
$F_{4}$ & 52 & $9$ & $(1,0,0,0)$ & $(0,0,0,1)$ & 2 & $\alpha_3,\alpha_4$\\
\hline
$G_{2}$ & 14 & $4$ & $(1,0)$ & $(0,1)$ & 3 & $\alpha_2$ \\
\hline
\end{tabular}
\caption{Basic Lie algebra data}
\label{tab:lie-algebra-data}
\end{table}

For completeness, we also give the form of the Cartan matrices $A$ in each of the cases.
By a slight abuse of notation, we indicate the Cartan matrix of an affine Lie algebra by the
label of the algebra itself.

The Cartan matrices for $A_r$ and $B_r$ read
\begin{align}
    A_r &= \begin{pmatrix}
        2 & -1 & 0 & \cdots & 0 & 0\\
        -1 & 2 & -1 & \cdots & 0 & 0\\
        0 & -1 & 2 & \cdots & 0 & 0\\
        \vdots & \vdots & \vdots & \ddots & \vdots & \vdots \\
        0 & 0 & 0 & \cdots & 2 & -1\\
        0 & 0 & 0 & \cdots & -1 & 2\\        
    \end{pmatrix}
    &
    B_r &= \begin{pmatrix}
    2 & -1 & 0 & \cdots & 0 & 0\\
    -1 & 2 & -1 & \cdots & 0 & 0\\
    0 & -1 & 2 & \cdots & 0 & 0\\
    \vdots & \vdots & \vdots & \ddots & \vdots & \vdots \\
    0 & 0 & 0 & \cdots & 2 & -2\\
    0 & 0 & 0 & \cdots & -1 & 2\\        
    \end{pmatrix} \ .
\end{align}
The Cartan matrices for $C_r$ and $D_r$ read
\begin{align}
    C_r &= \begin{pmatrix}
        2 & -1 & 0 & \cdots & 0 & 0\\
        -1 & 2 & -1 & \cdots & 0 & 0\\
        0 & -1 & 2 & \cdots & 0 & 0\\
        \vdots & \vdots & \vdots & \ddots & \vdots & \vdots \\
        0 & 0 & 0 & \cdots & 2 & -1\\
        0 & 0 & 0 & \cdots & -2 & 2\\        
    \end{pmatrix}
    &
    D_r &= \begin{pmatrix}
        2 & -1 & 0 & \cdots & 0 & 0 & 0\\
        -1 & 2 & -1 & \cdots & 0 & 0 & 0\\
        0 & -1 & 2 & \cdots & 0 & 0 & 0\\
        \vdots & \vdots & \vdots & \ddots & \vdots & \vdots \\
        0 & 0 & 0 & \cdots & 2 & -1 & -1\\
        0 & 0 & 0 & \cdots & -1 & 2 & 0\\        
        0 & 0 & 0 & \cdots & -1 & 0 & 2\\        
    \end{pmatrix} \ .
\end{align}
The Cartan matrices for $E_6$ and $E_7$ read
\begin{align}
    E_6 &= \begin{pmatrix}
        2 & -1 & 0 & 0 & 0 & 0\\
        -1 & 2 & -1 & 0 & 0 & 0\\
        0 & -1 & 2 & -1 & 0 & -1\\
        0 & 0 & -1 & 2 & -1 & 0\\
        0 & 0 & 0 & -1 & 2 & 0\\
        0 & 0 & -1 & 0 & 0 & 2\\
    \end{pmatrix}
&    
    E_7 &= \begin{pmatrix}
        2 & -1 & 0 & 0 & 0 & 0 & 0\\
        -1 & 2 & -1 & 0 & 0 & 0 & 0\\
        0 & -1 & 2 & -1 & 0 & 0 & -1\\
        0 & 0 & -1 & 2 & -1 & 0 & 0\\
        0 & 0 & 0 & -1 & 2 & -1 & 0\\
        0 & 0 & 0 & 0 & -1 & 2 & 0\\
        0 & 0 & -1 & 0 & 0 & 0 & 2\\
    \end{pmatrix} \ .
\end{align}
The Cartan matrices for $E_8$, $F_4$ and $G_2$ read
\begin{align}
    E_8 &= \begin{pmatrix}
        2 & -1 & 0 & 0 & 0 & 0 & 0 & 0\\
        -1 & 2 & -1 & 0 & 0 & 0 & 0 & 0\\
        0 & -1 & 2 & -1 & 0 & 0 & 0 & 0\\
        0 & 0 & -1 & 2 & -1 & 0 & 0 & 0\\
        0 & 0 & 0 & -1 & 2 & -1 & 0 & -1\\
        0 & 0 & 0 & 0 & -1 & 2 & -1 & 0\\
        0 & 0 & 0 & 0 & 0 & -1 & 2 & 0\\
        0 & 0 & 0 & 0 & -1 & 0 & 0 & 2\\
    \end{pmatrix}
&
    F_4 &= \begin{pmatrix}
        2 & -1 & 0 & 0\\
        -1 & 2 & -2 & 0\\
        0 & -1 & 2 & -1\\
        0 & 0 & -1 & 2\\
    \end{pmatrix}
&
    G_2 &= \begin{pmatrix}
        2 & -3\\
        -1 & 2\\
    \end{pmatrix} \ .
\end{align}

\section{Cases checked: unitarity}

In this appendix, we list the cases that we checked explicitly when
determining which RTCs are (pseudo-)unitary.
We use the rank $r$ and the level $k$ to label the uniform cases, while for the non-uniform cases, we use the rank $r$ and the label $\ell$.
The cases that were checked for unitarity are listed in table
\ref{tab:checked-unitarity}.

\begin{table}[ht]
\begin{tabular}[t]{| c | c | c | c |}
\hline
type & $r$ & $k$ & \# cases per rank\\
\hline
\multirow{2}{*}{$A_r$} & $1\leq r \leq 10$ & $1 \leq k \leq 20$ & $20$ \\ 
& $11\leq r \leq 20$ & $k=1,2$ & $2$ \\ 
\hline
\multirow{2}{*}{$B_r$} & $3\leq r \leq10$ & $1 \leq k \leq 20$ & $20$ \\ 
& $11 \leq r \leq 20$ & $k=1,2,3$ & $3$ \\ 
\hline
\multirow{2}{*}{$C_r$} & $2\leq r \leq 10$ & $1 \leq k \leq 20$ & $20$ \\ 
& $11\leq r \leq 20$ & $k=1,2,3$ & $3$ \\ 
\hline
\multirow{2}{*}{$D_r$} & $4\leq r \leq 10$ & $1 \leq k \leq 20$ & $20$ \\ 
& $11 \leq r \leq 20$ & $k=1,2,3$ & $3$ \\ 
\hline
$E_{r}$ & $r=6,7,8$ & $1 \leq k \leq 40$ & $40$ \\ 
\hline
$F_{4}$ & $r=4$ & $1 \leq k \leq 40$ & $40$ \\ 
\hline
$G_2$ & $r=2$ & $1 \leq k \leq 40$ & $40$ \\
\hline
\end{tabular}
\begin{tabular}[t]{| c | c | c | c |}
\hline
type & $r$ & $\ell$ & \# cases per rank\\
\hline
\multirow{2}{*}{$B_r$} & $3\leq r \leq 10$ & $\ell = 2r+1,2r+3,\ldots, 2r+39$ & $20$ \\ 
& $11 \leq r \leq 20$ & $\ell = 2r+1, 2r+3, 2r+5$ & $3$ \\
\hline
\multirow{2}{*}{$C_r$} & $2\leq r \leq 10$ & $\ell = 2r+1,2r+3,\ldots, 2r+39$ & $20$ \\ 
& $11\leq r \leq 20$ & $\ell = 2r+1, 2r+3, 2r+5$ & $3$ \\
\hline
$F_4$ & $r=4$ & $\ell = 13, 15,\ldots, 91$ & $40$\\
\hline
$G_2$ & $r=2$ & $\ell = 7, 8, 10, \ldots, 60, 61$ & $40$\\
\hline
\end{tabular}
\caption{The cases that were checked for (pseudo-)unitarity. Left panel: uniform cases, right panel: non-uniform cases.}
\label{tab:checked-unitarity}
\end{table}

\section{Cases checked: modularity}

In this appendix, we list the cases that we checked explicitly when
determining which RTCs are modular.
We use the rank $r$ and the level $k$ to label the uniform cases, while for the non-uniform cases, we use the rank $r$ and the label $\ell$.
The cases that were checked for modularity are listed in table
\ref{tab:checked-modularity}.

Because of stability issues when checking the modularity, we had to use arbitrary precision calculations, using 30 digits precision.
This lead to the fact that we could check fewer cases for modularity in comparison to unitarity. 
Because the modularity follows more regular patterns in comparison to the unitarity, we nevertheless have full confidence in our results
as stated in Sec.~\ref{sec:results}.

\begin{table}[ht]
\begin{tabular}[t]{| c | c | c | c |}
\hline
type & $r$ & $k$ & \# cases per rank\\
\hline
\multirow{2}{*}{$A_r$} & $1\leq r \leq 8$ & $1 \leq k \leq 20$ & $20$ \\ 
& $r=9,10$ & $1 \leq k \leq 10$ & $10$ \\ 
\hline
\multirow{2}{*}{$B_r$} & $3\leq r \leq10$ & $1 \leq k \leq 20$ & $20$ \\ 
& $11 \leq r \leq 20$ & $k=1,2$ & $2$ \\ 
\hline
\multirow{2}{*}{$C_r$} & $2\leq r \leq 8$ & $1 \leq k \leq 20$ & $20$ \\ 
& $r=9,10$ & $1 \leq k \leq 10$ & $10$ \\ 
\hline
\multirow{2}{*}{$D_r$} & $4\leq r \leq 10$ & $1 \leq k \leq 20$ & $20$ \\ 
& $11 \leq r \leq 20$ & $k=1,2$ & $2$ \\ 
\hline
$E_{r}$ & $r=6,7,8$ & $1 \leq k \leq 20$ & $20$ \\ 
\hline
$F_{4}$ & $r=4$ & $1 \leq k \leq 20$ & $20$ \\ 
\hline
$G_2$ & $r=2$ & $1 \leq k \leq 20$ & $20$ \\
\hline
\end{tabular}
\begin{tabular}[t]{| c | c | c | c |}
\hline
type & $r$ & $\ell$ & \# cases per rank\\
\hline
\multirow{4}{*}{$B_r$} & $3\leq r \leq 7$ & $\ell = 2r+1,2r+3,\ldots, 2r+39$ & $20$ \\ 
& $r=8$ & $\ell = 2r+1,2r+3,\ldots, 2r+29$ & $15$ \\ 
& $r=9,10$ & $\ell = 2r+1,2r+3,\ldots, 2r+19$ & $10$ \\
& $11 \leq r \leq 15$ & $\ell = 2r+1,2r+3, 2r+5$ & $3$ \\
\hline
\multirow{3}{*}{$C_r$} & $2\leq r \leq 7$ & $\ell = 2r+1,2r+3,\ldots, 2r+39$ & $20$ \\ 
& $r=8$ & $\ell = 2r+1,2r+3,\ldots, 2r+29$ & $15$ \\ 
& $r=9,10$ & $\ell = 2r+1,2r+3,\ldots, 2r+19$ & $10$ \\
\hline
$F_4$ & $r=4$ & $\ell = 13, 15,\ldots, 51$ & $20$\\
\hline
$G_2$ & $r=2$ & $\ell = 7, 8, 10, \ldots, 34, 35$ & $20$\\
\hline
\end{tabular}
\caption{The cases that were checked for modularity. Left panel: uniform cases, right panel: non-uniform cases.}
\label{tab:checked-modularity}
\end{table}

\end{document}